\documentclass[12pt,reqno]{amsart}
\usepackage{amsmath}
\usepackage{amssymb}
\usepackage{amsthm}
\usepackage[english]{babel}
\usepackage{amsmath}
\usepackage{amssymb}
\usepackage{amsthm}

\newtheorem{theorem}{Theorem}
\newtheorem{lemma}{Lemma}
\title[Estimates of stability by the number of summands]{Estimates of stability with respect to the number of summands for distributions of successive sums of independent
identically distributed vectors}

\def\le{\leqslant}
\def\ge{\geqslant}

\author[A.Yu. Zaitsev]{Andrei Yu. Zaitsev}
\address{St.~Petersburg Department of Steklov Mathematical Institute
\newline\indent
Fontanka 27, St.~Petersburg 191023, Russia\newline\indent
and St.Petersburg State University, 7/9 Universitetskaya nab., St. Petersburg,
199034 Russia}

\email{zaitsev@pdmi.ras.ru}

\thanks {This work was supported by the St. Petersburg International
Leonhard Euler Mathematical Institute, grant agreement No.
075-15-2022-289 dated 06.04.2022.}


\begin{document}
\begin{abstract}
Let $X_1,\dots, X_n,\dots$ be i.i.d.\
$d$-dimensional random vectors with common distribution $F$. Then $S_n =
X_1+\dots+X_n$ has distribution $F^n$ (degree is understood in the sense of
convolution). Let
$$
\rho_{\mathcal{C}_d}(F,G) = \sup_A |F\{A\} - G\{A\}|,
$$
where the supremum is taken over all convex subsets of $\mathbb R^d$. Basic
 result is as follows. For any nontrivial distribution $F$
there is $c_1(F)$ such that
$$
\rho_{\mathcal{C}_d}(F^n, F^{n+1})\leq \frac{c_1(F)}{\sqrt n}
$$
for any natural $n$. The distribution $F$ is called trivial if
it is concentrated on a hyperplane that does not contain the origin.
Clearly, for such $F$
$$
\rho_{\mathcal{C}_d}(F^n, F^{n+1}) = 1.
$$
A similar result for the Prokhorov distance is also obtained.
For any $d$-dimensional distribution~$F$ there is a  $c_2(F)>0$ that depends only on $F$ and such that
\begin{multline}\nonumber
(F^n)\{A\}\le (F^{n+1})\{A^{c_2(F)}\}+\frac{c_2(F)}{\sqrt{n}}\\ \text{and}\quad (F^{n+1})\{A\}\leq (F^n)\{A^{c_2(F)}\}+\frac{c_2(F)} {\sqrt{n}}
\end{multline}
for any
Borel set~$ A $ for all positive integers~$n$. Here $A^{\varepsilon }$ is $ \varepsilon $-neighborhood of the set~$ A $.
\end{abstract}
\maketitle

\section {Proximity of distributions of successive sums on convex sets}

Let $X_1$, $X_2$, \ldots, $X_n$, \ldots be independent identically distributed (i.i.d.) random vectors in the space $
{\mathbf R}^d $ with distribution $F$. Products and powers of measures will be understood in the sense of convolution: \;${GH=G*H}$, \,$H^m=H^{m*}$, \,$H^0=E=E_0$,
where $E_x$ is the distribution concentrated at point $x\in
{\mathbf R}^d $. Then $S_n =X_1+\dots+X_n$ has distribution~$F^n$.  We will study how different the distribution  $F^{n+1}$ is from the distribution  $F^n$, i.e., how much the distribution of the sum $S_n$ may be changed after adding another independent term to it.  It will be shown that the difference between these distributions is small, and it does not simply tend to zero as ${n\to\infty}$, but has order $O(n^{-1/2})$, standard for estimates in limit theorems of probability theory.

Sums of independent random variables and vectors is a classical object of probability theory. Since the study of the binomial distribution, which appeared in Bernoulli's scheme back in the eighteenth century, the properties of the distributions of sums of i.i.d.\ terms have been one of the main subjects of research. Under some (sometimes very restrictive) conditions, all possible limit distributions were found for the distributions of centered and normalized sums (see~\cite{AG, GK, P87}). In the results of the present paper the stability in the number of terms of the distributions of sums of i.i.d.\  vectors has been established
for arbitrary distributions of terms in finite-dimensional Euclidean spaces. Moreover, the obtained estimates have the optimal order $O(n^{-1/2})$ when comparing the values of the probabilities of getting into an arbitrary convex set for the sums of $n$ and $n+1$ summands. Surprisingly, in such a simple and natural formulation the problem was previously considered only in the author's publications (including joint ones), starting from the 1980s (see \cite{2, z80, z81, z81a, z83, z88}).

Such a statement of the problem naturally arises when considering the problem of Kolmogorov~\cite{23} about
estimating the accuracy of infinitely divisible approximation of distributions of sums of i.i.d.\ random variables.
Le~Cam~\cite{LC} showed that a natural infinitely divisible approximation for $F^n$ can be the accompanying compound Poisson distribution
$$
e(nF)=e^{-n}\sum_{s=0}^\infty
 \frac{n^sF^s}{s!},
$$
proposed by Gnedenko (see~\cite{GK, ip73}). It is clear that when estimating the closeness of the distributions $F^n$ and $e(nF)$ it is useful to be able to estimate the proximity of the distributions $F^n$ and $F^s$.

  Let us first introduce some notation. Let $\frak F_d$ denote the set
of probability distributions defined on the Borel $\sigma$-field $\mathcal{B}_d $ of
subsets of the Euclidean space~${\mathbf R}^d$.

Let's define the distances between the distributions
 $$\rho_{\mathcal{C}_d}(F, G) = \sup\limits_{A\in\mathcal{C}_d } {|F\{A\} - G\{A\}|},$$
$$\rho_{\rm TV}(F, G) = \sup\limits_{A\in\mathcal{B}_d } {|F\{A\} - G\{A\}|},$$
 where $\mathcal{C}_d $ is a collection of convex subsets, and $\mathcal{B}_d $ is a collection of Borel subsets of $
{\mathbf R}^d $.
In the one-dimensional case we use the notation $\rho(F, G) = \sup\limits_{x\in \mathbf{R}} {|F(x) - G(x)|}$ for the Kolmogorov distance (the uniform distance between cumulative distribution functions $F(\,\cdot\,)$ and $G(\,\cdot\,)$). It is clear that
$$
\rho(F, G)\le\rho_{\mathcal{C}_1}(F, G),\quad \rho_{\mathcal{C}_1}(F, G)\le2\,\rho(F, G).
$$
By the symbols $c$ and $c(\,\cdot\,)$ we generally denote various positive absolute constants and quantities that depend only on the argument in brackets. Distribution of random
vector $\xi $ will be denoted by $\mathcal L(\xi )$.

The following theorem is the first main result of this article.

\begin{theorem}\label{main} For any nontrivial distribution~$F$ there is a quantity $c(F)$ that depends only on $F$ and such that
\begin{equation}\label{e0} \rho_{\mathcal{C}_d}(F^n, F^{n+1}) \le\frac{c(F)}{\sqrt{n}} \end{equation} for all natural~$n$.
\end{theorem}

We call a distribution $F$  {\it trivial} if it is concentrated on an affine hyperplane that does not contain the origin. It is clear that for such $F$
\begin{equation}\label{e9} \rho_{\mathcal{C}_d}(F^n, F^{n+1}) =1. \end{equation} Triviality means that \eqref{e9} is trivially satisfied for trivial~$F$ since hyperplanes are convex sets and distributions $F^n$  and $F^{n+1}$ are concentrated on different disjoint hyperplanes.
In the one-dimensional case, trivial  distributions are distributions~$E_a$ concentrated at points $a\ne0$.

Theorem \ref{main} is a very general result. Inequalities \eqref{e0} and \eqref{e9} give a complete information on the proximity of the distributions $F^n$ and $F^{n+1}$ on arbitrary convex sets for arbitrary  distributions $F\in\mathfrak F_d$. The constant $c(F)$ in  inequality \eqref{e0} can be as large as you like if the distribution of $F$ is close to some trivial distribution.

In the one-dimensional case, the statement of Theorem \ref{main} is contained in \cite[Theorem 4.2 of Chapter V]{2}.
It is  known for non-degenerate Gaussian distributions $\Phi\in\mathfrak F_d$, and the estimate is valid even for distance in variation:
\begin{equation}\label{eq1}
  \rho_{\rm TV}(\Phi^n, \Phi^{n+1}) \le\frac{c(\Phi)}{\sqrt{n}}.
\end{equation}

This inequality may be derived using the following Lemma \ref{1} (see \cite{SZ}, \cite[Lemma 8]{PS}, as well as \cite[inequalities (1.3), (1.7)]{GNU}).

\begin{lemma} \label{1}
 Let $\Phi_k\in\mathfrak F_d$, $k=1,2$, be Gaussian distributions with nonsingular covariance matrices $\Sigma_k$ and means $b_k$.
   Then
   $$ \rho_{\rm TV}(\Phi_1, \Phi_2) \le\frac{1}{2}\Big(\bigl\|\Sigma_1^{-1/2}\Sigma_2\Sigma_1^{-1 /2}-I_d\bigr\|_F
   +\bigl\|\Sigma_2^{-1/2}(b_1-b_2)\bigr\|\Big), $$
where $\left\|\,\cdot\,\right\|_F$ is the Frobenius norm, and $I_d$ is the $d$-dimensional identity matrix.
\end{lemma}

In order to prove \eqref{eq1} one should apply Lemma \ref{1} with $\Sigma_1=n\Sigma$, $\Sigma_2=(n+1)\Sigma$, $b_1=nb$, $b_2=(n+1)b$, where $\Sigma$ and $b$ are covariance matrix and mean of the random vector $\xi$ with $\mathcal{L}(\xi)=\Phi$.

The monograph \cite{2} also contains other estimates of the proximity of $n$ and ${(n+1)}$-fold convolutions of one-dimensional distributions, including those with constants independent of the distribution~$F$. At the end of this section we will formulate some of these results. In recent joint works \cite{GZ18, GZ21} most of the mentioned results were extended to the values of distributions in Hilbert space on convex polyhedra, see also~\cite{z88}. The constants depend only on the number of half-spaces involved in the definition of the polyhedron.

Theorem \ref{main} will be deduced in a relatively elementary way using the following Lemma \ref{saz}, due to  Sazonov \cite{S68}, see also \cite{BR}, \cite{S81}.

\begin{lemma} \label{saz}
Let $F\in\mathfrak F_d$ be a  probability distribution with
   $$
   \int_{\mathbf{R}^{d}}\|x\|^3\,F\{dx\}<\infty,
   $$
   and $\Phi$ is a Gaussian distribution with the same covariance matrix and the same mean as the  distribution $F$.
Then there is a quantity $c(F)$ that depends only on $F$ and such that
$$ \rho_{\mathcal{C}_d}(F^n, \Phi^{n}) \le\frac{c(F)}{\sqrt{n}} $$ for all positive integers~$n$.
\end{lemma}

The binomial distribution with parameters $n,p$ may be represented in  the form
$$
B_{n,p}=\sum_{k=0}^{n}b_k(n,p)\,E_k,
$$
where$$
b_k(n, p)=C_n^k (1-p)^{n-k}p^k,
$$
and $C_n^k=\frac{n!}{k!(n-k)!}$ are binomial coefficients.
Let $\eta_{n,p}$ be a random variable with distribution $B_{n,p}$. It is well known that $\mathbf{E}\,\eta_{n,p}=np$, $\mathrm{Var}\,\eta_{n,p}=np(1-p)$.
We need the following lemma on the proximity of binomial distributions in variation.
\begin{lemma}\label{l3}For $0<p<1$ and any positive integer $n$, the following inequality holds:
\begin{equation}\label{e8}
\rho_{\rm TV}(B_{n,p}, B_{n+1,p})\le\frac{c(p)}{\sqrt{n}}.
\end{equation}
\end{lemma}

This lemma can also be considered as an estimate of the proximity of $n$ and $(n+1)$-fold convolutions, since $B_{n,p}=(B_{1,p})^n$.

\noindent{\bf Proof.} Ratio
\begin{equation}\label{e5}
\frac{b_k(n+1,p)}{b_k(n,p)}=\frac{n+1}{n+1-k}\cdot(1-p)
\end{equation}
increases monotonically as~$k$ increases.
At some ~$k$ there is a transition from ratio values not exceeding one to ratio values greater than one.
 Therefore, the difference between the distribution functions $B_{n,p}(x)- B_{n+1,p}(x)$ with increasing $x>0$ first increases from zero to the maximum value, and then decreases to zero. From the above it follows that
\begin{equation}\label{e6}
\rho_{\rm TV}(B_{n,p}, B_{n+1,p})=\rho(B_{n,p}, B_{n+1,p}).
\end{equation}
It is clear that $$B_{n+1,p}=B_{n,p}B_{1,p}=B_{n,p}((1-p)E_0+pE_1)=(1-p)B_{ n,p}+pE_1B_{n,p}.$$ Therefore,
\begin{equation}\label{e7}
\rho(B_{n,p}, B_{n+1,p})=p\,\rho(B_{n,p}, E_1B_{n,p})=p\,\max_k\mathbf{P}\{\eta_{n,p}=k\}\le\frac{c(p)}{\sqrt{n}}.
\end{equation}The last inequality in  \eqref{e7} is easily derived using the Stirling formula.
Inequality \eqref{e8} follows from \eqref{e6} and \eqref{e7}. $\square$
\medskip

We need the following property of the  distance $\rho_{\mathcal{C}_d}$.

\begin{lemma}[see \cite{Zol}]\label{l179} Let $ F, G, H\in\mathfrak F_d$ be arbitrary distributions. Then
$\rho_{\mathcal{C}_d}(FH, GH)\le \rho_{\mathcal{C}_d}(F,G)$.
\end{lemma}

\noindent{\bf Proof of Theorem \ref{main}.}
Without loss of generality, we can assume that the distribution  $F$ is non-trivial and is not concentrated on some proper subspace of~$
{\mathbf R}^d $. The proof actually uses induction on the dimension  $d$, taking into account the fact that if the distribution of $F$ is concentrated on some proper subspace and is trivial on it, then it is trivial on the space ${\mathbf R}^d$ itself.
It is easy to understand that there is $p$ such that $0<p<1$ and
\begin{equation}\label{88}F=(1-p)U+pV,
\end{equation}
where $U$ is a probability distribution with bounded support and a non-singular covariance matrix, and $V$ is some probability distribution. It is clear that the distributions  $U$ and $V$ can be chosen in such a way that $(1-p)U$ is the restriction of the measure $F$ to a centered ball of sufficiently large radius, and $pV$ to the complement to this ball.
The value $c(F)$ from the formulation of Theorem~\ref{main} will depend on $p$ and on the moments of the distribution $U$ up to the third order inclusive.
There are representations $$
F^n=\sum_{k=0}^{n}b_k(n,p)\,V^kU^{n-k},\quad F^{n+1}=\sum_{k=0}^{ n+1}b_k(n+1, p)\,V^kU^{n+1-k}.
$$

Introduce the distributions
$$
G_n=\sum_{k=0}^{n}b_k(n,p)\,V^kU^{n+1-k}.
$$
Let $\Phi $ be the Gaussian distribution with the same mean and covariance matrix as
distribution~$U$.

Applying Lemmas \ref{saz} and \ref{l179}, and inequality \eqref{eq1}, we obtain that, for $k\in\mathbf{Z}$, $0\le k< n$,
\begin{multline}
\label{e1}
\rho_{\mathcal{C}_d}(V^k U^{n-k},V^k U^{n+1-k})\le\rho_{\mathcal{C}_d}(U^{n-k}, U^{n+1-k})\\
\le\rho_{\mathcal{C}_d}(U^{n-k}, \Phi^{n-k})+\rho_{\mathcal{C}_d}(\Phi^{n-k}, \Phi^{n+1-k})+\rho_{\mathcal{C}_d}(\Phi^{n+1-k}, U^{n+1-k})\\
\le\frac{c(F)}{\sqrt{n-k}}+\frac{c(\Phi)}{\sqrt{n-k}}+\frac{c(F)}{\sqrt{n+1-k}}\le\frac{c(F)}{\sqrt{n-k}}.
\end{multline}
Therefore,
\begin{multline}
\label{est}
\rho_{\mathcal{C}_d}(F^n, G_n) \le\sum_{k=0}^{n}b_k(n,p)\,\rho_{\mathcal{C}_d}(V^k U^{n-k},V^k U^{n+1-k})\\
\le b_n(n,p)+\sum_{k=0}^{n-1}b_k(n,p)\,\frac{c(F)}{\sqrt{n-k}}\\ \le b_n(n,p)+ c(F)\,\mathbf{E}\,\frac{\mathbf{1}\{\eta_{n,p}<n\}}{\sqrt{n-\eta_{n,p}}}.
\end{multline}
Here $\mathbf{1}\big\{A\big\}$ is the indicator of an event $A$.

According to Bernstein's inequality (see \cite[Theorem 4.1 of Chapter~I]{2}),
\begin{equation}
\label{bern}\mathbf{P}\big\{\eta_{n,p}-np\ge np(1-p)\big\}\le\exp(-np(1-p)/4).
\end{equation}
It is easy to see that $0<p(2-p)<1$ for $0<p<1$. Hence,
\begin{multline}
\label{bound}
\mathbf{E}\,\frac{\mathbf{1}\{\eta_{n,p}<n\}}{\sqrt{n-\eta_{n,p}}}=\mathbf{E}\,\frac{1}{\sqrt{n-\eta_{n,p}}}\,\mathbf{1}\big\{\eta_{n,p}<np+np(1-p)\big\}\\
+\mathbf{E}\,\frac{\mathbf{1}\{\eta_{n,p}<n\}}{\sqrt{n-\eta_{n,p}}}\,\mathbf{1}\big\{\eta_{n,p}\ge np+np(1-p)\big\}\\
\le\frac{c(p)}{\sqrt{n}}+\exp(-np(1-p)/4)\le\frac{c(p)}{\sqrt{n}}=\frac{c(F)}{\sqrt{n}}.
\end{multline}
 It is clear that
\begin{equation}\label{e46}
\rho_{\mathcal{C}_d}(G_n, F^{n+1})\le\sum_{k=0}^{n+1}\left|b_k(n,p)-b_k(n+1,p)\right|=2\,\rho_{\rm TV}(B_{n,p}, B_{n+1,p})
\end{equation}
(of course, we assume $b_{n+1}(n,p)=0$). Moreover,  $b_{n}(n,p)=p^n\le c(p)/\sqrt{n}$. It remains to apply \eqref{est},  \eqref{bound}, \eqref{e46} and Lemma \ref{l3}. $\square$
\medskip

A point $a$ is called the $q$-quantile of a one-dimensional distribution~$F$ if $F\{(-\infty,a)\} \le q$ and $F\{( a, \infty)\}\le 1 - q$, where $0 \le q \le 1$. For $q=1/2$, the $q$-quantile is called the median of   distribution $F$.
Let the point
$0$ be the $q$-quantile of the distribution $F$.
Then the following estimate for the Kolmogorov distance is valid (see \cite[Theorem 4.1 of Chapter V]{2}, as well as \cite{z80}):
\begin{equation}
\label{estimation}
\rho(F^n, F^{n+1}) \le \frac{c}{\sqrt{n\min{\{q, 1-q\}}}} \le \frac{c}{\sqrt{nq(1-q)}}.
\end{equation}
The dependence on $n$ and $q$ in this inequality is correct, since a similar lower estimate is valid, that is, the estimate \eqref{estimation} is optimal (see \cite[Example 4.1 of Chapter V]{2}).
This estimate is based on a special case of the Kolmogorov--Rogozin inequality (see \cite[Theorem 2.4 of Chapter~II]{2}).
In  inequality $(\ref{estimation})$ the absolute constant $c$ can be taken equal to $c_0=\frac{1+2\sqrt{2\pi}}{e^{3/8}}\approx 4.132847$ \cite{g18}.

If $0$ is the median of the distribution $F$, then
\begin{equation}
\label{estimation1}
\rho(F^n, F^{n+1}) \le \frac{c_0\sqrt2}{\sqrt{n}}=\frac{c}{\sqrt{n}}
\end{equation}
(see \cite{z81}). In particular, this is true if the distribution of $F$ is symmetric.
For symmetric  $F$, the following unexpected and paradoxical inequality is also true:
\begin{equation}
\label{estimation2}
\rho(F^n, F^{n+2}) \le \frac{c}{{n}}
\end{equation}
(see \cite[Theorem 5.2 of Chapter V]{2}). For  standard normal distribution~$F$ it follows from Lemma \ref{1}.
If $m$ is the median of the distribution~$F$, then the distribution $FE_{-m}$ has  zero median and
\begin{equation}
\label{estimation2}
\rho((FE_{-m})^n, (FE_{-m})^{n+1})=\rho(F^n, F^{n+1}E_{-m})  \le \frac{c}{\sqrt{n}}.
\end{equation}
It is clear that if $F=E_a$ is a degenerate distribution with $a\ne0$, then $$\rho(F^n, F^{n+1})=1.$$ In any other case $\rho (F^n, F^{n+1})\le c(F)/\sqrt n$.
Indeed,
\begin{multline}
\label{estimation3}
\rho(F^n, F^{n+1})\le\rho(F^n, F^{n+1}E_{-m}) +\rho(F^{n+1}, F^{n+1}E_{-m})\\
 \le \frac{c}{\sqrt{n}}+Q(F^{n+1}, |m|)\le c(F)/\sqrt n.
\end{multline}
Here $Q(\,\cdot\,,\,\cdot\,)$ is the L\'evy concentration function, and the last inequality follows from the Kolmogorov--Rogozin inequality.
Thus, we  obtained a one-dimensional version of Theorem \ref{main}.

For $d=1$, Lemma \ref{saz} follows from the well-known Berry--Esseen inequality. It is clear that for distributions with finite moments of the third order the statement of Theorem~\ref{main} can be easily deduced from \eqref{eq1} and Lemma~\ref{saz} using the triangle inequality.

\section{Estimates of the Prokhorov distance}

In this section we  formulate an analogue of Theorem \ref{main} for the Prokho\-rov distance
 \cite{hh19} metrizing the weak convergence of probability distributions (see Theorem \ref{proh0} below). The question about the possibility of obtaining such an analogue was raised by Youri  Davydov during the author's talk concerning Theorem \ref{main}.

 The Prokhorov distance between distributions $G,H\in\mathfrak F_d$ is defined as
\begin{multline}
\pi (G,H)=\inf \Big\{ \varepsilon>0 :G\{A\}\le H\{A^{\varepsilon}\}+\varepsilon, H\{A\}\leq G\{A^{\varepsilon}\}+\varepsilon\\
 \text{for any
Borel set~$ A $}\Big\},\nonumber
\end{multline}
where
$A^{\varepsilon }=\{y\in \mathbf{R}^{d}:\inf\limits_{x\in A}\left\|
x-y\right\| <\varepsilon \}$ denotes the $ \varepsilon $-neighborhood of a set~$ A\in\mathcal{B}_d $.

We need the following property of the Prokhorov distance.

\begin{lemma}[see \cite{Zol}]\label{l177} Let $ F, G, H\in\mathfrak F_d$ be arbitrary distributions. Then
$\pi(FH, GH)\le \pi(F,G)$.
\end{lemma}

The following Lemma~\ref{l1} is usually called the Strassen--Dudley theorem (see \cite{hh9, Schay, Str}).

\begin{lemma}\label{l1} Let $ F, G\in\mathfrak F_d$ be arbitrary distributions. Then
\begin{multline}\pi(F, G)=\inf \Big\{ \varepsilon>0 :\hbox{ one can construct}\\ \hbox{on the same probability space}\\
\hbox{ the random vectors $\xi$ and $\eta$ with
$\mathcal{L}(\xi)=F$ and $\mathcal{L}(\eta)=G$} \\ \hbox{ so that }
\mathbf{P}\left\{ \Vert \xi -\eta \Vert >\varepsilon \right\}\le \varepsilon
 \Big\}.\end{multline}
\end{lemma}

If $X$ is a random vector with distribution $F$ and $a>0$, we will denote by $F_{(a)}$ the distribution of the normalized random vector $X/\sqrt{a}$.
The following Lemma~\ref{l12} can be easily derived using Lemma~\ref{l1}.

\begin{lemma}\label{l12} Let $ F, G\in\mathfrak F_d$ be arbitrary distributions. Then, for  any $a,b>0$,
$$\pi(F_{(b)},G_{(b)})\le\max\Big\{\frac{\sqrt{a}}{\sqrt{b}}, 1\Big\}\,\pi(F_{(a)},G_{(a)}).$$
\end{lemma}

The proof  of Theorem \ref{proh0} uses induction on the dimension $d$.
Other steps in the proof almost literally repeat the proof of Theorem~\ref{main} in Section~1. Only instead of Lemma~\ref{saz} we should use the following
Lemma~\ref{l4}, due to V.V.~Yurinskii \cite{yu75}.

\begin{lemma} \label{l4} Let $F\in\mathfrak F_d$ be a  probability distribution with
   $$
   \int_{\mathbf{R}^{d}}\|x\|^3\,F\{dx\}<\infty,
   $$
   and $\Phi$ is a Gaussian distribution with the same covariance matrix and the same mean as the  distbution $F$. Then
     there exists a quantity $c(F)$ that depends only on $F$ and such that
$$ \pi(F_{(n)}^n, \Phi_{(n)}^{n}) \le\frac{c(F)}{\sqrt{n}} $$ for all natural numbers~$ n$.
\end{lemma}
The value $c(F)$ from the formulation of Lemma~\ref{l4} is depending on the moments of the distribution $F$ up to the third order inclusive.
The original formulation of Lemma \ref{l4} of  Yurinskii \cite{yu75} is a little bit different. The random vectors are normalised not only by $\sqrt{n}$ but also by $\sigma$, where $\sigma^2$ is the maximal eigenvalue of the covariance matrix of summands. In order to obtain the statement of Lemma \ref{l4} one should use in addition Lemma~\ref{l12}.

 The second main
 result is    Theorem \ref{proh0}. We will prove this theorem in Section~3. Theorem \ref{proh} is an auxiliary result concerning the most part of non-degenerate distributions~$F$.

 \begin{theorem}\label{proh} Assume that $F\in\mathfrak F_d$ is a probability distribution such that
\begin{equation}\label{88}F=(1-p)U+pV,
\end{equation}
where $0<p<1$, $U\in\mathfrak F_d$ is a probability distribution with bounded support and a non-singular covariance matrix, and $V\in\mathfrak F_d$ is some probability distribution. Then there exists a quantity $c(F)$ that depends only on $F$ and such that
\begin{equation}\label{e01} \pi(F_{(n)}^n, F_{(n)}^{n+1}) \le\frac{c(F)}{\sqrt{n }} \end{equation} for all natural~$n$.
\end{theorem}

The value  $c(F)$ from the statement of  Theorem \ref{proh} depends on $p$ and  on the moments of the distribution $U$ up to the third order inclusive.

\begin{theorem}\label{proh0} For any distribution $F\in\mathfrak F_d$ there exists a quantity $c(F)$ that depends only on $F$ and such that
\begin{equation}\label{e015} \pi(F_{(n)}^n, F_{(n)}^{n+1}) \le\frac{c(F)}{\sqrt{n }} \end{equation} for all natural~$n$.
\end{theorem}

The right-hand side of inequality~\eqref{e015} has the correct order in~$n$. To verify this, it is enough to take as $F$ the symmetric one-dimensional distribution $F=E_{-1}/2+E_{1}/2$ and even $n$. Then the distributions $F_{(n)}^n$ and $F_{(n)}^{n+1}$ are concentrated, respectively, on the sets
$\big\{2k/\sqrt{n}, k\in \mathbf{Z}\big\}$ and $\big\{(2k+1)/\sqrt{n}, k\in \mathbf{ Z}\big\}$, and
  \begin{equation}\label{e77}
\max_k\mathbf{P}\{S_{n}=k\}=\mathbf{P}\{S_{n}=0\}\ge\frac{c}{\sqrt{n}}.
\end{equation}The last inequality in the formula \eqref{e77} is easily derived using the Stirling formula. From the above it follows that
\begin{equation}\label{e0155} \pi(F_{(n)}^n, F_{(n)}^{n+1}) \ge\frac{c}{\sqrt{n }}. \end{equation}

In Theorems \ref{proh} and \ref{proh0}  we do not divide distributions into trivial and non-trivial.
Note that for all $G,H\in\mathfrak F_d$
\begin{equation}\label{e0145}
\pi(G, H) \le \rho_{\rm TV}(G, H) .
\end{equation}
Therefore, for non-degenerate Gaussian distributions, Theorem \ref{proh0} follows from  inequality~\eqref{eq1}.

It is easy to see that, for all $a, b\in\mathbf{R}^d$,
\begin{equation}\label{e014}
\pi(E_a, E_b) =\min\big\{1,\|b-a\|\big\}.
\end{equation}

It is clear that $F_{(n)}^n=\mathcal{L}(S_n/\sqrt{n})$ is not close to the degenerate distribution~$E_0$ unless the distribution~$F$ is degenerate.  Indeed, using the Kolmogorov--Rogozin inequality for the L\'evy concentration functions (see  \cite[Theorem 2.4 of Chapter II]{2}), one can  show that, for any non-degenerate distribution $F$, there exists a $c(F)$ such that $$\mathbf{P}\big\{\|S_n/\sqrt{n}\|\le c(F)\big\}\le1/2.$$ At the same time it is obvious that one can choose normalizing constants $\varphi(n)$ to be so large that  $\pi\big(\mathcal{L}(S_n/\varphi(n)),
\mathcal{L}(S_{n+1}/\varphi(n))\big)$ is small due to the fact that both distributions are close to the degenerate distribution $E_0$.

Thus, the normalization by $\sqrt{n}$ is natural when considering the Pro\-kho\-rov distance, which is not invariant under scale transformation.
For distributions of non-normalized sums, the statement of Theorem \ref{proh0} may be not true in general. In particular, $\pi(F^n, F^{n+1})=1$ for $F=E_a$ with $\|a\|\ge 1$ (see \eqref{e014}).
\medskip

From Theorem \ref{proh0}, Lemma~\ref{l1} and from the definition of the Prokhorov distance, the following Theorems \ref{proh2} and \ref{proh3} may be easily deduced.

\begin{theorem}\label{proh2} For any distribution $F\in\mathfrak F_d$ there exists a quantity $c(F)$  depending only on $F$ and such that
\begin{multline}
(F^n)\{A\}\le (F^{n+1})\{A^{c(F)}\}+\frac{c(F)}{\sqrt{n}}\\ \text{and}\quad (F^{n+1})\{A\}\leq (F^n)\{A^{c(F)}\}+\frac{c(F)} {\sqrt{n}}
\end{multline}
for any
Borel set~$ A $ and for all positive integers~$n$.
\end{theorem}

\begin{theorem}\label{proh3} For any distribution $F\in\mathfrak F_d$ there exists a quantity $c(F)$  depending only on $F$ and such that,
for any natural~$n$,
one can construct on the same probabilistic
space random vectors \,$\xi_n $ \,and \,$\eta_n $ \,with
\,$\mathcal{L}(\xi_n )=F^{n+1}$ \,and \,$\mathcal{L}(\eta_n )=F^n$, \,so that
\begin{equation}
\mathbf{P}\left\{ \Vert \xi_n -\eta_n \Vert >c(F) \right\} \le \frac{c(F)}{\sqrt{n}}.  \label{pi1}
\end{equation}
\end{theorem}

Note that the vectors $\xi_n=S_{n+1}$ and $\eta_n=S_n$ have the required distributions, but for the vector $\xi_n -\eta_n=X_{n+1}$ in this case  inequality~\eqref{pi1} does not hold, of course,  if the distribution $F$ has unbounded support. If the support of  distribution $F$ is bounded, then the statement of Theorem~\ref{proh3} for $\xi_n=S_{n+1}$ and $\eta_n=S_n$ is obvious, and the right-hand side of  inequality \eqref{pi1} may be replaced by zero.

Theorems \ref{proh0}--\ref{proh3} are also very general statements. They describe the closeness of the distributions $F^n$ and~$F^{n+1}$ on arbitrary Borel sets for arbitrary distributions $F\in\mathfrak F_d$. In fact, the statements of Theorems \ref{proh0}, \ref{proh2} and \ref{proh3} are equivalent. Note that Theorems \ref{proh2} and \ref{proh3} say about the closeness of the distributions $F^{n+1}$ and ~$F^n$ of non-normalized vectors $S_{n+1}$ and $S_n $, which once again indicates the naturality of choosing normalization by $\sqrt{n}$ considering the distributions $F_{(n)}^{n}$ and ~$F_{(n)}^{n+1}$ of vectors $S_{n}/\sqrt{n}$ and $S_{n+1}/\sqrt{n}$ in the formulation of  Theorem \ref{proh0}.

Let $X_1, X_2, \ldots$ be i.i.d.\ random vectors with a common distribution $F\in\mathfrak F_d$ and let $\mu$ be an integer valued non-negative random variable independent of the sequence $\{X_j\}_{j=1}^\infty$. Denote
$
  G=\mathcal{L}(X_1+\cdots+X_\mu)
$.
It is well known that then
\begin{equation}
G=\sum_{k=0}^\infty\mathbf{P}\{\mu=k\}\,F^k.
   \label{008a}\end{equation}
It is clear that estimates of the proximity of distributions $F^{n+1}$ and~$F^n$ can be useful when comparing distributions of the form~\eqref{008a} (see \cite[\S 5 of Chapter V]{ 2}, \cite{GZ18, GZ21}).

\section{Proofs}

\noindent{\bf Proof of Theorem \ref{proh}.}
We have representations
\begin{equation}\label{8899}F_{(n)}=(1-p)U_{(n)}+pV_{(n)},
\end{equation}
 $$
F_{(n)}^n=\sum_{k=0}^{n}b_k(n,p)\,V^k_{(n)}U_{(n)}^{n-k},\quad F_{(n)}^{n+1}=\sum_{k=0}^{ n+1}b_k(n+1, p)\,V_{(n)}^kU_{(n)}^{n+1-k}.
$$

Introduce the distributions
$$
G_n=\sum_{k=0}^{n}b_k(n,p)\,V_{(n)}^kU_{(n)}^{n+1-k}.
$$
Let $\Phi $ be the Gaussian distribution with the same mean and covariance matrix as
distribution~$U$. Applying inequalities \eqref{eq1}, \eqref{e0145} and Lemmas \ref{l177},~\ref{l12},~\ref{l4}, we obtain that, for $k\in\mathbf{Z}$, $0\le k< np(2-p)$,
\begin{multline}
\label{e1}
\pi(V_{(n)}^k U_{(n)}^{n-k},V_{(n)}^k U_{(n)}^{n+1-k})\le\pi(U_{(n)}^{n-k}, U_{(n)}^{n+1-k})\\
\le\pi(U_{(n)}^{n-k}, \Phi_{(n)}^{n-k})+\pi(\Phi_{(n)}^{n-k}, \Phi_{(n)}^{n+1-k})+\pi(\Phi_{(n)}^{n+1-k}, U_{(n)}^{n+1-k})\\
\le\pi(U_{(n-k)}^{n-k}, \Phi_{(n-k)}^{n-k})+\pi(\Phi_{(n-k)}^{n-k}, \Phi_{(n-k)}^{n+1-k})+c\,\pi(\Phi_{(n+1-k)}^{n+1-k}, U_{(n+1-k)}^{n+1-k})\\
\le\frac{c(F)}{\sqrt{n-k}}+\frac{c(\Phi)}{\sqrt{n-k}}+\frac{c(F)}{\sqrt{n+1-k}}\le\frac{c(F)}{\sqrt{n-k}}<\frac{c(F)}{\sqrt{n}}:=\varepsilon_n.
\end{multline}
We used again  that $0<p(2-p)<1$ for $0<p<1$.

Recall that  $\eta_{n,p}$ is a random variable with binomial distribution $B_{n,p}$. It satisfies inequality~\eqref{bern}.
Using \eqref{bern}, \eqref{e1}, we get, for any Borel set~$A$,
\begin{multline}
\label{est}
F_{(n)}^n\{A\}=\sum_{k=0}^{n}b_k(n,p)\,V_{(n)}^kU_{(n)}^{n-k}\{A\} \\ =\sum_{\substack{k=0\\ k\ge np(2-p)}}^{n}b_k(n,p)\,V_{(n)}^kU_{(n)}^{n-k}\{A\}
+ \sum_{\substack{k=0\\ k< np(2-p)}}^{n}b_k(n,p)\,V_{(n)}^kU_{(n)}^{n-k}\{A\}\\
\le\mathbf{P}\big\{\eta_{n,p}\ge np(2-p)\big\} + \sum_{\substack{k=0\\ k< np(2-p)}}^{n}b_k(n,p)\,\big(V_{(n)}^kU_{(n)}^{n+1-k}\{A^{\varepsilon_n}\}+\varepsilon_n\big)\\
\le\exp(-np(1-p)/4)+G_n\{A^{\varepsilon_n}\}+\varepsilon_n\\
\le  G_n\{A^{\varepsilon_n}\}+\frac{c(F)}{\sqrt{n}}.\nonumber
\end{multline}
Similarly,
$$
G_n\{A\}\le F_{(n)}^n\{A^{\varepsilon_n}\}+\frac{c(F)}{\sqrt{n}}.
$$
Hence,
\begin{equation}\label{e47}
\pi(F_{(n)}^{n},G_n)\le\frac{c(F)}{\sqrt{n}}.
\end{equation}

It is clear that
\begin{equation}\label{e4}
\pi(G_n, F_{(n)}^{n+1})\le\sum_{k=0}^{n+1}\left|b_k(n,p)-b_k(n+1,p)\right|=2\,\rho_{\rm TV}(B_{n,p}, B_{n+1,p})
\end{equation}
(of course, we assume $b_{n+1}(n,p)=0$). It remains to apply \eqref{e47},  \eqref{e4} and Lemma~\ref{l3}.  $\square$\medskip

\noindent{\bf Proof of Theorem \ref{proh0}.} We will use induction on the dimension  $d$. Let $d=1$ and $F=E_a$ for $a\in\mathbf{R}$. Then $F_{(n)}^n=E_{na/\sqrt{n}}$,
$F_{(n)}^{n+1}=E_{(n+1)a/\sqrt{n}}$ and, using \eqref{e014}, we have $$\pi(F_{(n)}^n, F_{(n)}^{n+1}) =\pi(E_0, E_{a/\sqrt{n}})\le|a|/\sqrt{n},$$  proving  \eqref{e015} in this case. If $F$ is non-degenerate, $F\ne E_a$, we can represent $F$ in the form \eqref{88} in such a way that $(1-p)U$ is the restriction of the measure $F$ to an interval $[-T,T]$, and $pV$ to $\mathbf{R}\setminus[-T,T]$, the complement to this interval. Choosing $T$ large enough, we can guarantee that a random variable $\xi$ with $\mathcal{L}(\xi)=U$ has non-zero variance. Now one-dimensional version of Theorem~\ref{proh0} follows from Theorem~\ref{proh}.

Assume that Theorem~\ref{proh0} is proved for $(d-1)$-dimensional distributions~$F$. Let us prove it in the $d$-dimensional case.
It is easy to understand that  we can represent $F$ in the form \eqref{88} with some $p$ such that $0<p<1$ and
where $U\in\mathfrak F_d$ is a probability distribution with bounded support and a non-singular covariance matrix, and $V\in\mathfrak F_d$ is some probability distribution. It is clear that the distributions  $U$ and $V$ can be chosen in such a way that $(1-p)U$ is the restriction of the measure $F$ to a centered ball of sufficiently large radius, and $pV$ to the complement to this ball. Moreover, if for all radii the distribution $U$ has some singular covariance matrix, then distribution $F$ is concentrated on an affine  hyperplane~$H$. If this hyperplane contains the origin, then it is a $(d-1)$-dimensional linear subspace of $\mathbf{R}^d$ and the statement of theorem follows from the induction hypothesis.
If $H$ does not  contain the origin, then it is a $(d-1)$-dimensional linear subspace $H_0$ shifted by a vector $a\in\mathbf{R}^d$: $H=H_0+a$.
In this case, the distribution $F$ may be represented as $F=GE_a$, where the distribution $G$ is concentrated on the  hyperplane~$H_0$. By the induction hypothesis, \begin{equation}\label{88999}\pi(G_{(n)}^n, G_{(n)}^{n+1}) \le\frac{c(G)}{\sqrt{n }}.\end{equation}
Using \eqref{e014}, \eqref{88999} and Lemma~\ref{l177}, we get
\begin{multline}\label{mm}
 \pi(F_{(n)}^n, F_{(n)}^{n+1}) =\pi(G_{(n)}^nE_{na/\sqrt{n}},G_{(n)}^{n+1} E_{(n+1)a/\sqrt{n}}) \\
  \le \pi(G_{(n)}^nE_{na/\sqrt{n}},G_{(n)}^{n+1} E_{na/\sqrt{n}}) +\pi(E_{na/\sqrt{n}}, E_{(n+1)a/\sqrt{n}})
   \\
  \le \frac{c(G)}{\sqrt{n }}+\frac{\|a\|}{\sqrt{n }} \le \frac{c(F)}{\sqrt{n }}.\nonumber
\end{multline}
Theorem~\ref{proh0} is proved. $\square$

The author is grateful to Yu. A. Davydov for his question about the Prokhorov distance and to V. V. Ulyanov for useful advices on bibliographic issues.

\end{document}